\documentclass[a4paper]{article}
\usepackage{hyperref}
\usepackage{amsmath}
\usepackage[dvips]{graphicx}
\usepackage{amsthm}
\usepackage{epsfig}
\usepackage{amssymb}
\usepackage{dsfont}
\usepackage{float}

\newtheorem{theorem}{Theorem}[section]

\newtheorem{pro}{Proposition}[section]

\newtheorem{de}{Definition}[section]

\begin{document}

\begin{center}
{\Large   Existence of peakons for a cubic generalization of the Camassa-Holm equation}
\end{center}
\vskip 5mm

\begin{center}
{\sc Yun Wang and Lixin Tian} \\

\vskip2mm
Institute of Mathematics,
School of Mathematical Sciences,\\
Nanjing Normal University,
Nanjing, 210023, PR China
\vskip 5mm
\end{center}

\vskip 5mm {\leftskip5mm\rightskip5mm \normalsize
\noindent{\bf{Abstract}} In this paper, we study the following generalized Camassa-Holm equation with both 
cubic and quadratic nonlinearities:
$$
m_{t}+k_{1}(3uu_{x}m+u^2m_{x})+k_{2}(2mu_{x}+m_{x}u)=0, \quad m=u-u_{xx},
$$
which is presented as a linear combination of the Novikov equation and the Camassa-Holm equation with constants
$k_{1}$ and $k_{2}$. The model is a cubic generalization of the Camassa-Holm equation. It is shown that 
the equation  admits single-peaked soliton and periodic peakons.
}

\newtheorem{proposition}[theorem]{Proposition}

\renewcommand{\theequation}{\thesection.\arabic{equation}}
\catcode`@=11
\@addtoreset{equation}{section}
\catcode`@=12

\section{Introduction}

The well-known Camassa-Holm(CH) equation
\begin{equation}\label{CH}
m_{t}+um_{x}+2u_{x}m=0,\quad m=u-u_{xx},
\end{equation}
which was proposed by Camassa and Holm as a nonlinear model for the unidirectional propagation of 
the shallow water waves over a flat bottom with $u(x, t)$ representing the water's free surface \cite{CL, CH,RS}. 
It has attracted much attention in the past decades.
In addition, the CH equation (\ref{CH}) has several nice geometrical structures due to, for example, its description about a 
geodesic flow on the diffeomorphism group on the circle \cite{SK} and its derivation from a non-stretching invariant planar 
curve flow in the centro-equiaffine geometry \cite{CQ}. Moreover, well-posedness theory and wave breaking phenomenon 
of the CH equation were studied extensively, and many interesting results have been deduced, see 
\cite{AC, ACJ, ACJ1, ACJ2, LP1, XZ1, XZ2}. The stability and interaction of peakons were discussed in several references
\cite{ACW, ACW1, FGLQ}. Among these properties, a remarkable one is that it
admits the single peakons and periodic peakons in the following forms
\begin{equation}\label{SP}
    \varphi_{c}(x,t)=ce^{-|x-ct|}, \quad c\in \mathbb{ R},
\end{equation}
and
\begin{equation}\label{PP}
    u_{c}(x,t)=\frac{c}{\operatorname{sh}(1/2)}\operatorname{ch}(\frac{1}{2}-(x-ct)+[x-ct]), 
    \quad c\in \mathbb{ R},
\end{equation}
where the notation $[x]$ denotes the largest integer part of the real number $x\in\mathbb{R}$.

In addition to the CH equation being an integrable model with peakons, other integrable peakon models, which include the 
Degasperis-Procesi equation and the cubic nonlinear peakon equations \cite{AHJ, Nv}, have been found. Indeed, two integrable 
CH-type equations with cubic nonlinearity have been discovered recently. The first one is mCH equation:
\begin{equation}\label{mCH}
m_{t}+\left((u^2-u^2_{x})m\right)_{x}=0, ~~m=u-u_{xx},
\end{equation}
and the second one is the so-called  Novikov equation:
\begin{equation}\label{Nv}
u_{t}-u_{txx}+4u^2u_{x}=3uu_{x}u_{xx}+u^2u_{xxx}, \quad t>0,  \quad x\in\mathbb{R}.
\end{equation}
The perturbative symmetry approach \cite{MN} yielded necessary conditions for PDEs to admit 
infinitely many symmetries. Using this approach, Novikov \cite{Nv} was able to isolate Eq.(\ref{Nv}) in a symmetry classification 
and also found its first few symmetries.
He subsequently found a scalar Lax pair for it, and also proved that the equation is integrable.
Hone and Wang \cite{AHJ} showed that 
equation (\ref{Nv}) arised as a zero curvature equation $F_{t}-G_{x}+[F, G]=0$,  which is the compatibility condition for the 
linear system
$$
\left \{
   \begin{array}{l}
   \Psi_{x}=F\Psi,\\
   \Psi_{t}=G\Psi,
   \end{array}
   \right.
$$
where $m=u-u_{xx}$,
$$
F=
\left(
\begin{array}{ccc}
0&m\lambda&1\\
0&0&m\lambda\\
1&0&0\\
\end{array}
\right),
G=
\left(
\begin{array}{ccc}
\frac{1}{3\lambda^2}-uu_{x}&\frac{u_{x}}{\lambda}-u^2m\lambda &u^2_{x}\\
\frac{u}{\lambda}&-\frac{2}{3\lambda^2}&-\frac{u}{\lambda}-u^2m\lambda\\
-u^2&\frac{u}{\lambda}&\frac{1}{3\lambda^2}+uu_{x}\\
\end{array}
\right).
$$
They also proved that it was related to a negative flow in the Sawada-Kotera hierarchy by a reciprocal transformation.
 By defining a new dependent 
variable $m$, Eq.(\ref{Nv}) can be written as 
\begin{equation}\label{Nv1}
m_{t}+u^2m_{x}+3uu_{x}m=0,\quad m=u-u_{xx}.
\end{equation}

Analogous to the Camassa-Holm equation, the Novikov equation has a bi-Hamiltonian structure and an infinite sequence 
of conserved quantities. In addition, the single peakon for the Novikov equation was obtained in \cite{AHJ}, which takes the form
$$
u(t,x)=\pm\sqrt{c}e^{-|x-ct|}, \quad c>0,
$$
and the periodic peakons \cite{WT}
$$
 u_{c}(x,t)=\sqrt{c}\frac{\operatorname{ch}(\frac{1}{2}-(x-ct)+[x-ct])}{\operatorname{ch}(1/2)}, 
    \quad c>0.
$$
Then, Liu, Liu and Qu \cite{LL1} proved the single peakons are orbital stable.  Wang, Tian \cite{WT} also proved the existence and orbital stability 
of periodic peakons.


On the other hand, applying tri-Hamiltonian duality to the modified Korteweg-de Vries (mKdV) equation leads to the 
modified Camassa-Holm (mCH) equation with cubic nonlinearity.
More generally, applying tri-Hamiltonian duality to the bi-Hamiltonian Gardner equation
\begin{equation}{G}
u_{t}+u_{xxx}+k_{1}u^2u_{x}+k_{2}uu_{x}=0,
\end{equation}
the resulting dual system is the following generalized modified Camassa-Holm (gmCH) equation with both cubic and quadratic
nonlinearities \cite{BF}:
\begin{equation}\label{Ll}
 m_{t}+k_{1}\left((u^2-u^2_{x})m\right)_{x}+k_{2}(2u_{x}m+um_{x})=0,~~~m=u-u_{xx}.
 \end{equation}
Recently, it was found \cite{QXL} that, for $k_{1}\neq0$,  the gmCH equation (\ref{Ll}) admits a single peakon of 
the form
$$
\varphi_{c}(t, x)=ae^{-|x-ct|}, ~~c\in \mathbb{R},
$$
with
$$
\displaystyle a=\frac{3}{4}\frac{-k_{2}\pm\sqrt{k^2_{2}+\frac{8}{3}k_{1}c}}{k_{1}},~~k^2_{2}+\frac{8}{3}k_{1}c\geq0,
$$
and also found \cite{LL} that, for $k_{1}\neq0$,  the gmCH equation (\ref{Ll}) admits periodic peakons of 
the form
$$
u_{c}(t, x)=a \operatorname{ch}(\frac{1}{2}-(x-ct)+[x-ct]), 
$$
where
$$
\displaystyle a=\frac{3}{4}\frac{-k_{2}\operatorname{ch}(1/2)\pm\sqrt{k^2_{2}\operatorname{ch}^2(1/2)+\frac{4}{3}k_{1}
c(1+2\operatorname{ch}^2(1/2))}}{k_{1}(1+2\operatorname{ch}^2(1/2))}
$$
and 
$$
k^2_{2}\operatorname{ch}^2(1/2)+\frac{4}{3}k_{1}c(1+2\operatorname{ch}^2(1/2))\geqslant0.
$$

 The existence of (periodic) peakons is of interest for the nonlinear integrable equations since they are relatively
 new solitary waves. More importantly, in the theory of water waves many papers have investigated the 
 Stokes waves of greatest height, traveling waves which are smooth everywhere except at the crest where the lateral
 tangents differ. There is no closed form available for these waves, and the peakons capture the essential features of the 
 extreme waves-waves of great amplitude that are exact solutions of the governing equations for irrotational water waves,
 see the discussion in \cite{A, AJ}.
Inspired by \cite{LL, QXL}, we focus on the following generalized Camassa-Holm equation with both 
cubic and quadratic nonlinearities:
\begin{equation}\label{gCH}
m_{t}+k_{1}(3uu_{x}m+u^2m_{x})+k_{2}(2mu_{x}+m_{x}u)=0, \quad m=u-u_{xx},
\end{equation}
where $k_{1}$ and $k_{2}$ are arbitrary constants. It is clear that equation (\ref{gCH}) reduces to the CH 
equation for $k_{1}=0, k_{2}=1$ and the Novikov equation for $k_{1}=1, k_{2}=0$,  respectively. 
Equation (\ref{gCH}) is actually a linear combination of CH equation (\ref{CH}) and cubic nonlinear 
equation (\ref{Nv1}). Therefore, we may view equation (\ref{gCH}) as a generalization of the CH equation, or 
simply call equation  (\ref{gCH}) a generalized CH equation.
Like the Camassa-Holm and Novikov equations, the new equation admits peaked 
soliton solutions. 

\section{Preliminaries}
In this paper, we are concerned with the Cauchy problem for the generalized CH equation on both line and the unit circle:
\begin{equation}\label{2.1}
\left \{
   \begin{array}{l}
   m_{t}+k_{1}(3uu_{x}m+u^2m_{x})+k_{2}(2mu_{x}+m_{x}u)=0, \quad t>0,~ x\in X=\mathbb{R}~ or~ \mathcal{S},\\
   m(t, x)=u(t, x)-u_{xx}(t,x),\\
   u(0, x)=u_{0}(x), \quad x\in X.
  \end{array}
   \right.
\end{equation}
First, we will require the notion of strong (or classical) solutions as follows:
\begin{de}\label{de2.1}
If $u\in C\left([0, T); H^s(X)\right)\cap C^1([0, T); H^{s-1}(X))$ with $s>\frac{5}{2}$ and some $T>0$ 
satisfies (\ref{2.1}),  then $u$ is  called a strong solution on $[0,T)$.  If $u$ is a strong solution on $[0,T)$ for 
every $T>0$, then it is called a global strong solution.
\end{de}

The following local well-posedness result and properties for strong solutions on the line and unit circle can be
established using the same approach as in \cite{GLQ}. 
\begin{pro}
Let $u_{0}\in H^s(X)$ with $s>\frac{5}{2}$. Then there exists a time $T>0$ such that the initial value problem 
(\ref{2.1}) has a unique strong solution $u\in C\left([0, T); H^s(X)\right)\cap C^1([0, T); H^{s-1}(X))$ and the map
$u_{0}\rightarrow u$ is continuous from a neighborhood of $u_{0}$ in $H^s(X)$ into $u\in C\left([0, T); H^s(X)\right)
\cap C^1([0, T); H^{s-1}(X))$.
\end{pro}

If $m=u-u_{xx}$ is substituted in terms of $u$ into the generalized CH equation (\ref{2.1}), then the resulting
fully nonlinear partial differential equation takes the following form:
\begin{equation}\label{2.2}
\begin{array}{ll}
\displaystyle u_{t}+k_{1}u^2u_{x}+\frac{1}{2}k_{1}(1-\partial^2_{x})^{-1}u^3_{x}+k_{1}(1-\partial^2_{x})^{-1}\partial_{x}(u^3
+\frac{3}{2}uu^2_{x})\\[2mm]
\qquad \qquad\qquad\qquad\displaystyle+k_{2}uu_{x}+k_{2}\partial_{x}(1-\partial^2_{x})^{-1}(u^2+\frac{1}{2}u^2_{x})=0.
\end{array}
\end{equation}
Taking the convolution with the Green's function for the Helmholtz operator $(1-\partial^2_{x})$, equation (\ref{2.2})
can be rewritten as
\begin{equation}\label{p2.2}
\begin{array}{ll}
\displaystyle u_{t}+k_{1}u^2u_{x}+\frac{1}{2}k_{1}G(x)*u^3_{x}+k_{1}G(x)*\partial_{x}(u^3
+\frac{3}{2}uu^2_{x})\\[2mm]
\qquad \qquad\qquad\qquad\displaystyle+k_{2}uu_{x}+k_{2}G(x)*\partial_{x}(u^2+\frac{1}{2}u^2_{x})=0.
\end{array}
\end{equation}
Note that $u$ can be formulated by the Green function $G(x)$ as
\begin{equation}\label{2.3}
u=(1-\partial^2_{x})^{-1}m=G*m,
\end{equation}
where $G(x)=\frac{1}{2}e^{-|x|}$ for the non-periodic case, $G(x)=\displaystyle\frac{\operatorname{ch}(1/2-x+[x])}
{2\operatorname{sh}(1/2)}$ for the periodic case, and $*$ denotes the convolution product on $X$, defined
by
$$
(f*g)(x)=\displaystyle\int_{X}f(y)g(x-y)dy.
$$

The above formulation (\ref{2.2}) allows us to define the periodic weak solutions as follows.
\begin{de}
Given initial data $u_{0}\in W^{1,3}(\mathcal{S})$, the function $u\in L^{\infty}_{loc}([0, T), W^{1,3}_{loc}
(\mathcal{S}))$ is called a periodic weak solution to the initial value problem (\ref{2.1}) if it satisfies the 
following identity:
\begin{equation}
\begin{array}{ll}
\displaystyle\int^T_{0}\int_{\mathcal{S}}\Big[u\partial_{t}\phi+\frac{k_{1}}{3}u^3\partial_{x}\phi+k_{1}G(x)*\left(u^3
+\frac{3}{2}uu^2_{x}\right)\partial_{x}\phi-k_{1}G(x)*\left(\frac{u^3_{x}}{2}\right)\phi\\[3mm]
\qquad+\displaystyle\frac{k_{2}}{2}u^2\partial_{x}\phi+k_{2}G(x)*(u^2+\frac{1}{2}u^2_{x})\partial_{x}\phi\Big]dxdt
+\displaystyle\int_{\mathcal{S}}
u_{0}(x)\phi(0,x)dx=0,
\end{array}
\end{equation}
for any smooth test function $\phi(t, x)\in C^{\infty}_{c}\left([0, T)\times \mathcal{S}\right)$. If $u$ is a weak 
solution on $[0, T)$ for every $T>0$, then it is called a global periodic weak solution.
\end{de}

\section{Peakon solutions}
In this section, we derive the single-peaked solutions and periodic peakons of equation  (\ref{gCH}).
\subsection{Single-peaked solutions}
Firstly, we give the result on the existence of single-peaked solutions for the generalized CH equation (\ref{gCH}).
\begin{theorem}(Single peakons)\label{s2.1}
For the wave speed $c$ satisfying $k^2_{2}+4k_{1}c\geq0$, equation (\ref{gCH}) with $k_{1}\neq0$ admits the single 
peakons of the form:
\begin{equation}\label{s2.5}
u=Ae^{-|x-ct|},
\end{equation}
where $\displaystyle A=\frac{-k_{2}\pm\sqrt{k^2_{2}+4k_{1}c}}{2k_{1}}$.
\end{theorem}
$Proof$.
Firstly,  let us suppose the single-peaked solutions of equation (\ref{gCH}) in the form of 
\begin{equation}\label{s3.1}
u=Ae^{-|x-ct|},
\end{equation}
where $A$ is to be determined.
The derivatives of expression (\ref{s3.1}) do not exist at $x=ct$, thus (\ref{s3.1}) can not satisfy equation (\ref{gCH}) in 
the classical sense. However, in the weak sense, we can write out the expressions of $u_{x}$ and $u_{t}$  under the help of 
distribution:
\begin{equation}\label{s3.2}
u_{x}=-A sgn(x-ct)e^{-|x-ct|},~~u_{t}=cA sgn(x-ct)e^{-|x-ct|}.
\end{equation}
Next, we need consider two cases (i) $x>ct$ and (ii) $x<ct$. \\
For $x>ct$, we calculate from (\ref{s3.1}) and (\ref{s3.2}) that
\begin{equation}\label{s3.3}
u_{t}+k_{1}u^2u_{x}+k_{2}uu_{x}=Ace^{-(x-ct)}-k_{1}A^3e^{-3(x-ct)}-k_{2}A^2e^{-2(x-ct)}.
\end{equation}
Note that the Green function $G(x)=\frac{1}{2}e^{-|x|}$ in the non-periodic case, it is thus deduced that 
\begin{equation}\label{s3.4}
\begin{array}{ll}
\displaystyle\frac{1}{2}k_{1}G(x)*u^3_{x}+k_{1}G(x)*\partial_{x}(u^3+\frac{3}{2}uu^2_{x})\\[3mm]
\quad=\displaystyle-4k_{1}A^3\left(-\int^{ct}_{-\infty}+\int^x_{ct}+\int^{+\infty}_{x}\right)e^{-|x-y|-3|y-ct|}dy\\[3mm]
\quad=\displaystyle-k_{1}A^3e^{-(x-ct)}+k_{1}A^3e^{-3(x-ct)},
\end{array}
\end{equation}
and 
\begin{equation}\label{s3.5}
\begin{array}{ll}
\displaystyle k_{2}G(x)*\partial_{x}(u^2+\frac{1}{2}u^2_{x})\\[3mm]
\quad=\displaystyle-\frac{3}{2}k_{2}A^2\int_{R}sgn(y-ct)e^{-|x-y|-2|y-ct|}dy\\[3mm]
\quad=\displaystyle-\frac{3}{2}k_{2}A^2\left(-\int^{ct}_{-\infty}+\int^x_{ct}+\int^{+\infty}_{x}\right)e^{-|x-y|-2|y-ct|}dy\\[3mm]
\quad=\displaystyle-k_{2}A^2e^{-(x-ct)}+k_{2}A^2e^{-2(x-ct)}.
\end{array}
\end{equation}
The case $x<ct$ is similar to $x>ct$,  here we do not compute in detail.\\
Plugging  (\ref{s3.3}), (\ref{s3.4}) and (\ref{s3.5}) into (\ref{p2.2}), we deduce that 
\begin{equation}\label{s3.6}
\begin{array}{ll}
\displaystyle u_{t}+k_{1}u^2u_{x}+\frac{1}{2}k_{1}G(x)*u^3_{x}+k_{1}G(x)*\partial_{x}(u^3
+\frac{3}{2}uu^2_{x})\\[2mm]
\qquad \qquad\qquad\qquad\displaystyle+k_{2}uu_{x}+k_{2}G(x)*\partial_{x}(u^2+\frac{1}{2}u^2_{x})\\[3mm]
\quad=(Ac-k_{1}A^3-k_{2}A^2)e^{-(x-ct)}\\[3mm]
\quad=0.
\end{array}
\end{equation}
Therefore, we are able to conclude from (\ref{s3.6}) that $A$ should satisfy
\begin{equation}
k_{1}A^2+k_{2}A-c=0.
\end{equation}
In general, we obtain
\begin{equation}
A=\frac{-k_{2}\pm\sqrt{k^2_{2}+4k_{1}c}}{2k_{1}}
\end{equation}
with $k^2_{2}+4k_{1}c\geq0$ and $k_{1}\neq0$.

\subsection{Periodic peakons}
The following theorem shows the existence of periodic peakons for the generalized CH equation (\ref{gCH}).
\begin{theorem}(Periodic peakons)\label{th2.1}
For the wave speed $c$ satisfying $k^2_{2}\operatorname{ch}^2(1/2)+4k_{1}c\left(1+\operatorname{sh}
^2(1/2)\right) \geq0$, equation (\ref{gCH}) with $k_{1}\neq0$ possesses the periodic 
peakons of the form:
\begin{equation}
u_{c}(x, t)=a\operatorname{ch}(\zeta),~~\zeta=\frac{1}{2}-(x-ct)+[x-ct],
\end{equation}
where
\begin{equation}\label{2.6}
a=\frac{-k_{2}\operatorname{ch}(1/2)\pm\sqrt{k^2_{2}\operatorname{ch}^2(1/2)+4k_{1}c\left(1+\operatorname{sh}
^2(1/2)\right)}}{2k_{1}\left(1+\operatorname{sh}^2(1/2)\right)}
\end{equation}
as the global periodic weak solutions to (\ref{2.1}) in the sense of Definition 2.2.
\end{theorem}

$Proof$.
Firstly, we identify $\mathcal{S}=[0, 1)$ and regard $u_{c}(t, x)$ as spatial periodic function on $\mathcal{S}$ 
with period one.  On one hand, it is noted that $u_{c}$ is continuous on $\mathcal{S}$ with peak at $x=0$. 
On the other hand, $u_{c}$ is smooth on $(0,1)$ and for all $t\in \mathbb{R}^+$,
\begin{equation}\label{3.1}
\partial_{x}u_{c}(t, x)=-a\operatorname{sh}(\zeta)\in L^{\infty}(\mathcal{S}).
\end{equation}
Hence, we denote $u_{c, 0}(x)=u_{c}(0, x)$ with $x\in\mathcal{S}$, then it holds that
\begin{equation}\label{3.2}
\lim_{t\rightarrow0^+}\|u_{c}(t, \cdot)-u_{c, 0}(\cdot)\|_{W^{1, \infty}(\mathcal{S})}=0.
\end{equation}
As in (\ref{3.1}), it is found that
\begin{equation}\label{3.3}
\partial_{t} u_{c}(x, t)=ac\operatorname{sh}(\zeta)\in L^\infty(\mathcal{S}),~t\geq0.
\end{equation}
Using (\ref{3.1})-(\ref{3.3}) and integration by parts, it is thus deduced that, for every test function $\phi(t, x)\in 
C^\infty_{c}\left([0, \infty)\times\mathcal{S} \right)$,
\begin{equation}\label{3.4}
\begin{array}{ll}
\displaystyle\int^\infty_{0}\int_{\mathcal{S}}\left(u_{c}\partial_{t}\phi+\frac{k_{1}}{3}u_{c}^3\partial_{x}\phi
+\frac{k_{2}}{2}u_{c}^2\partial_{x}\phi\right)dxdt+\int_{\mathcal{S}}u_{c,0}(x)\phi(x, 0)dx \\[4mm]
\quad\quad\quad\quad=-\displaystyle\int^\infty_{0}\int_{\mathcal{S}}\phi(\partial_{t}u_{c}+k_{1}u_{c}^2
\partial_{x}u_{c}+k_{2}u_{c}\partial_{x}u_{c})dxdt\\[4mm]
\quad\quad\quad\quad=\displaystyle\int^\infty_{0}\int_{\mathcal{S}}\phi \left((-ac
+k_{1}a^3)\operatorname{sh}(\zeta)+k_{1}a^3\operatorname{sh}^3(\zeta)
+k_{2}a^2\operatorname{sh}(\zeta)\operatorname{ch}(\zeta)\right)dxdt.
\end{array}
\end{equation}
On the other hand, noticing from the explicit form of the Green function $G(x)$ for the periodic case that
$$
G(x)=\frac{\operatorname{ch}(1/2-x+[x])}{2\operatorname{sh}(1/2)}~~ and~~
G_{x}(x)=-\frac{\operatorname{sh}(1/2-x+[x])}{2\operatorname{sh}(1/2)}, ~ x\in\mathbb{R},
$$
it follows from (\ref{3.1}), (\ref{3.3}) and the proof of Theorem 4.1 in \cite{Qu} that
\begin{equation}\label{3.5}
\begin{array}{ll}
\displaystyle\int^\infty_{0}\int_{\mathcal{S}}\left[ k_{1}G(x)*\left(u_{c}^3+\frac{3}{2}u_{c}(\partial_{x}
u_{c})^2\right)\partial_{x}\phi-\frac{k_{1}}{2}G(x)*(\partial_{x}u_{c})^3\phi\right]dxdt\\[4mm]
\quad\quad\quad\quad=\displaystyle k_{1}a^3\int^\infty_{0}\int_{\mathcal{S}}\phi G(x)*\left(3\operatorname{sh}(\zeta)
+\frac{7}{2}\operatorname{sh}
^3(\zeta)\right)dxdt\\[4mm]
\qquad\qquad\qquad\displaystyle-\frac{3}{2}k_{1}a^3\displaystyle\int^\infty_{0}\int_{\mathcal{S}}\phi G_{x}(x)*\left
(\operatorname{ch}(\zeta)\operatorname{sh}^2(\zeta)\right)dxdt\\[4mm]
\quad\quad\quad\quad=\displaystyle k_{1}a^3\int^\infty_{0}\int_{\mathcal{S}}\phi\left(\operatorname{sh}^2(1/2)\cdot
\operatorname{sh}(\zeta)-\operatorname{sh}^3(\zeta)\right)dxdt.
\end{array}
\end{equation}

Next, we compute directly that
\begin{equation}\label{3.23}
\begin{array}{ll}
\displaystyle\int^\infty_{0}\int_{\mathcal{S}}k_{2}G(x)*\left(u_{c}^2+\frac{1}{2}(\partial_{x}u_{c})^2\right)\partial_{x}
\phi dxdt\\[4mm]
\qquad\displaystyle=\frac{3k_{2}}{2}a^2\int^\infty_{0}\int_{\mathcal{S}}\phi G(x)*\operatorname{sh}(2\zeta)dxdt.
\end{array}
\end{equation}
To prove, we consider two cases: (i) $x>ct$ and (ii) $x<ct$. When $x>ct$, a direct calculation gives rise to
\begin{equation}\label{3.24}
\begin{array}{ll}
G(x)*\operatorname{sh}(2\zeta)(t, x)\\[4mm]
\quad\displaystyle=\frac{1}{2\operatorname{sh}(1/2)}\int_{\mathcal{S}}\operatorname{ch}(1/2-(x-y)+[x-y])\cdot
\operatorname{sh}\left(1-2(y-ct)+2[y-ct]\right)dy\\[4mm]
\quad\displaystyle=\frac{1}{2\operatorname{sh}(1/2)}\Big[\int^{ct}_{0}\operatorname{ch}(1/2-x+y)\cdot
\operatorname{sh}(-1-2y+2ct)dy\\[4mm]
\qquad\qquad\qquad\qquad\displaystyle+\int^{x}_{ct}\operatorname{ch}(1/2-x+y)\cdot\operatorname{sh}
(1-2y+2ct)dy\\[4mm]
~\qquad\qquad\qquad\qquad\displaystyle+\int^{1}_{x}\operatorname{ch}(1/2+x-y)\cdot\operatorname{sh}(1-2y+2ct)dy
 \Big]\\[4mm]
\quad\displaystyle=\frac{2}{3}\left[\operatorname{ch}(1/2)\operatorname{sh}(1/2-(x-ct))-
\operatorname{sh}(1/2-(x-ct)) \operatorname{ch}\left(1/2-(x-ct)\right)\right].
\end{array}
\end{equation}
In a similar manner, for $x<ct$,
\begin{equation}\label{3.25}
\begin{array}{ll}
G(x)*\operatorname{sh}(2\zeta)(t, x)\\[4mm]
\quad\displaystyle=\frac{2}{3}\left[-\operatorname{ch}(1/2)\operatorname{sh}(1/2+(x-ct))+\operatorname{sh}(1/2+(x-ct)) 
\operatorname{ch}\left(1/2+(x-ct)\right)\right].
\end{array}
\end{equation}
Plugging (\ref{3.24}) and (\ref{3.25}) into (\ref{3.23}), it is deduced by a straightforward computation that
\begin{equation}\label{3.26}
\begin{array}{ll}
\displaystyle\int^\infty_{0}\int_{\mathcal{S}}k_{2}G(x)*\left(u_{c}^2+\frac{1}{2}(\partial_{x}u_{c})^2\right)\partial_{x}
\phi dxdt\\[4mm]
\qquad\displaystyle=k_{2}a^2\int^\infty_{0}\int_{\mathcal{S}}\phi\left(\operatorname{ch}(1/2)\operatorname{sh}(\zeta)
-\operatorname{sh}(\zeta)\operatorname{ch}(\zeta) \right)dxdt.
\end{array}
\end{equation}

In view of (\ref{3.4}), (\ref{3.5}) and (\ref{3.26}), we have
\begin{equation}\label{3.27}
\begin{array}{ll}
\displaystyle\int^\infty_{0}\int_{\mathcal{S}}[u_{c}\partial_{t}\phi+\frac{k_{1}}{3}u_{c}^3\partial_{x}\phi+\frac{k_{2}}{2}
u_{c}^2\partial_{x}\phi\\[3mm]
\qquad\qquad\displaystyle+k_{1}G(x)*\left(u_{c}^3+\frac{3}{2}u_{c}(\partial_{x}u_{c})^2\right)\partial_{x}
\phi-k_{1}G(x)*\left(\frac{(\partial_{x}u_{c})^3}{2}\right)\phi\\[3mm]
\qquad\qquad\qquad\displaystyle+k_{2}G(x)* (u_{c}^2+\frac{1}{2}(\partial_{x}u_{c})^2)   
\partial_{x}\phi]dxdt+\displaystyle\int_{\mathcal{S}}u_{c, 0}(x)\phi(0,x)dx\\[3mm]
\quad\displaystyle=\int^\infty_{0}\int_{\mathcal{S}}\phi a\Big[k_{1}(1+\operatorname{sh}^2(1/2))a^2
+k_{2}\operatorname{ch}(1/2)a-c\Big]\operatorname{sh}(\zeta)dxdt.
\end{array}
\end{equation}

If a takes value as (\ref{2.6}), then
$$
k_{1}(1+\operatorname{sh}^2(1/2))a^2+k_{2}\operatorname{ch}(1/2)a-c=0,
$$
which implies that
\begin{equation}\label{3.28}
\begin{array}{ll}
\displaystyle\int^\infty_{0}\int_{\mathcal{S}}[u_{c}\partial_{t}\phi+\frac{k_{1}}{3}u_{c}^3\partial_{x}\phi+\frac{k_{2}}{2}
u_{c}^2\partial_{x}\phi\\[3mm]
\qquad\displaystyle+k_{1}G(x)*\left(u_{c}^3+\frac{3}{2}u_{c}(\partial_{x}u_{c})^2\right)\partial_{x}
\phi-\frac{k_{1}}{2}G(x)*(\partial_{x}u_{c})^3\phi\\[3mm]
\quad\qquad\displaystyle+k_{2}G(x)* (u_{c}^2+\frac{1}{2}(\partial_{x}u_{c})^2)   
\partial_{x}\phi]dxdt+\displaystyle\int_{\mathcal{S}}u_{c, 0}(x)\phi(0,x)dx=0,
\end{array}
\end{equation}
for any test function $\phi(x, t)\in C^\infty_{c}([0, \infty)\times \mathcal{S})$. Thus the theorem is proved.

\end{document}